\documentclass{amsart}
\usepackage{amssymb, amsmath, amsfonts, amscd}

\input xypic

\theoremstyle{plain}
\newtheorem{theorem}{Theorem}[section]

\newtheorem{lemma}[theorem]{Lemma}
\newtheorem{proposition}[theorem]{Proposition}

\theoremstyle{definition}
\newtheorem{definition}[theorem]{Definition}

\theoremstyle{remark}

\numberwithin{equation}{theorem}

\newcommand{\F}{\mathcal{F}}
\newcommand{\G}{\mathcal{G}}
\newcommand{\I}{\mathcal{I}}

\newcommand{\E}{\mathcal{E}}
\renewcommand{\O}{\mathcal{O} }

\renewcommand{\P}{\mathbf{P} }
\renewcommand{\Pr}{\mathcal{P} }
\newcommand{\Diff}{\mathcal{D} }
\renewcommand{\u}{\frac{1}{u} }
\renewcommand{\t}{\frac{1}{t} }

\newcommand{\mb}[1]{\mathbf{#1} }

\newcommand{\Spec}{\operatorname{Spec} }

\newcommand{\Proj}{\operatorname{Proj} }

\newcommand{\GL}{\operatorname{GL} }

\newcommand{\Z}{\mathbf{Z} }

\newcommand{\p}{\mathbf{P}}

\begin{document}

\title{Principal parts on the projective line over arbitrary rings}
\author{Helge Maakestad }
\address{Dept. of Mathematics, KTH, Stockholm}
\email{makestad@math.kth.se}
\thanks{}
\date{Spring 2004}
\begin{abstract} We develop techniques to split explicitly the sheaf of principal parts $\Pr^k(\O(n))$ as left and right
module on the projective line over an arbitrary ring.
We then apply the techniques developed to split
the principal parts $\Pr^k(\O(n))$ for all  $n\in \Z$ and $k\geq 1$ as left and right $\O$-module on the projective line 
over any field of characteristic zero giving a complete description of the principal parts on the projective line.
\end{abstract}

\maketitle

\tableofcontents
\section{Introduction}
In previous papers (see \cite{maa1},\cite{maa2}) the structure of the principal parts on the
projective line was studied using explicit and multilinear
algebra techniques, and a partial classification was obtained on the projective line and projective space. 
In this paper the techniques from \cite{maa1} are generalized to study the splitting type of the principal parts $\Pr^k(\O(n))$ 
for all $k\geq 1$ and all $n\in \Z$ as left and right module on the projective line over any ring. 
The techniques are powerful enough to enable us to give a complete classification of the splitting type of the principal parts
on the projective line over any field of characteristic zero. 
The main result of the paper is the following:
\begin{theorem}\label{projectiveline} Let $F$ be a field of characteristic zero. Consider the principal parts $\Pr^k(\O(d))$ 
on the projective line $\P^1_F$. The structure of the principal parts as left module is as follows:
If $k\geq 1$ and $n<0$ or $n\geq k$ there exist an isomorphism
\begin{equation}
\Pr^k(\O(n))^{left}\cong \oplus^{k+1} \O(n-k) 
\end{equation} 
of left $\O$-modules.
If $0\leq n <k$ there exist an isomorphism
\begin{equation}
\Pr^k(\O(n))^{left}\cong \O^{n+1}\oplus \O(-k-1)^{k-n}
\end{equation} 
of left $\O$-modules.
The structure of the principal parts as right module is as follows: 
For all $n\in \Z$ and $k\geq 1$ there exist an isomorphism
\begin{equation}
\Pr^k(\O(n))^{right}\cong \O(n)\oplus \O(n-k-1)^k
\end{equation} 
of right $\O$-modules. 
\end{theorem} 
\begin{proof} See Theorem \ref{leftspl1},\ref{leftspl2} and \ref{rightspl1}.
\end{proof} 
It is hoped that the techniques developed in the paper will shed light on the problem of classifying the 
structure of the principal parts as left and right module on the projective line  over any field. 

The paper is organized as follows: In section one we define and prove general properties of the sheaf of principal parts.
In section two we calculate the transition matrices of the principal parts $\Pr^k(\O(n))$ for all $k\geq 1$ and $n\in \Z$ on 
the projective line over the integers. In section three we develop techniques to split the principal parts as left and
right module over any ring. In section four we apply the techniques developed in section three to give a complete
classification of the splitting type of the principal parts on the projective line over any field of characteristic zero.

\section{Sheaves of principal parts}

In this section we define and give general properties of the sheaves of principal parts following \cite{egaIV}, chapter 16. 

\textbf{Notation}. Let in the following $X/S$ be a separated scheme, and consider
the two projection maps $p,q:X\times X\rightarrow X$. The diagonal embedding $\Delta :X\rightarrow X\times X$ 
is a closed immersion, and we get an exact sequence of sheaves of $\O_{X\times X}$-modules on $X\times X$
\begin{equation}  \label{prinexact} 
  0\rightarrow \I^{k+1} \rightarrow \O_{X\times X} \rightarrow \O_{\Delta^k} \rightarrow 0, 
\end{equation} 
where $\I$ is the sheaf of ideals defining the diagonal in $X\times X$.

\begin{definition} \label{principalparts} 
Let $\E$ be a quasi-coherent $\O_X$-module. We define
the \emph{$k$'th order sheaf of principal parts} of $\E$, to be
\[ \Pr^k_X(\E)=p_*(\O_{\Delta^k}\otimes_{X\times X} q^*\E).\]
We write $\Pr^k_X$ for the module $\Pr^k_X(\O_X)$.
\end{definition}

When it is clear from the context which scheme we are working on, we 
write $\Pr^k(\E)$ instead of $\Pr^k_X(\E)$. 

\begin{proposition} \label{fundamentalexact} Let $X/S$ be smooth and let $\E$ be a locally free
$\O_X$-module.
There exists an exact sequence
\[ 0 \rightarrow S^k(\Omega^1_X)\otimes \E \rightarrow \Pr^k_X(\E) \rightarrow
\Pr^{k-1}_X(\E) \rightarrow 0 \]
of left $\O_X$-modules, where $k=1,2,\dots$.
\end{proposition}
\begin{proof} See \cite{laksov}, section 4. 
\end{proof}
From Proposition \ref{fundamentalexact} it follows that for a smooth morphism $X\rightarrow S$ of relative dimension
$n$, and $\E$ a locally free sheaf on $X$ of rank $e$, the principal
parts sheaf $\Pr^k(\E)$ is locally free of rank $e\binom{n+k}{n}$.

\begin{proposition} \label{functoriality}
Let $f:X\rightarrow Y$ be a map of smooth schemes
over $S$, and let 
$\E$ be a locally free $\O_Y$-module. 
There exists a commutative diagram of exact sequences
\[ \diagram 0 \rto & S^k(f^*\Omega^1_Y)\otimes f^*\E \dto \rto
& f^*\Pr^k_Y(\E) \dto \rto & f^*\Pr^{k-1}_Y(\E) \dto \rto & 0 \\
0 \rto & S^k(\Omega^1_X)\otimes f^*\E \rto & \Pr^k_X(f^*\E) \rto & \Pr^{k-1}_X(f^*\E)
\rto & 0 \enddiagram \]
of left $\O_X$-modules for all $k=1,2,\dots$.
\end{proposition}
\begin{proof} See \cite{perkinson}. 
\end{proof}
Note that it follows from  Proposition \ref{functoriality} that for any open set $U\subseteq X$ of a scheme $X$ smooth 
over $S$ there is an isomorphism
\begin{equation}\label{local}
\Pr^k_X(\E)|_U\cong \Pr^k_U(\E|_U)
\end{equation}
of left $\O_U$-modules: Let $i:U\rightarrow X$ be the inclusion map. It follows that $i$ is smooth. By induction and applying
the 5-lemma Proposition \ref{functoriality} proves  the isomorphism \ref{local}.
Hence the sheaf of principal parts can be studied using local calculations.

Given any $\O_X$-modules $\F$ and $\G$ one may define the \emph{sheaf of differential operators} of order $k$
from $\F$ to $\G$ following \cite{egaIV} section 16.8, denoted $\Diff_X ^k(\F,\G)$, and there exists an isomorphism
\begin{equation} \label{prindiff}
  Hom _X(\Pr_X ^k( \F ),\G)\cong \Diff_X ^k(\F ,\G ) 
\end{equation} 
of sheaves of abelian groups. The sheaf of abelian groups $\Diff_X^k(\F,\G)$ is naturally a sheaf of  $\O_X$-bimodules and it follows 
from the isomorphism (\ref{prindiff}) that the sheaf of principal parts $\Pr_X^k(\F)$ is also a sheaf of $\O_X$-bimodules.
This means that for any open set $U$ of $X$, and local sections $a,b\in \O(U)$ there exists 
a left and right multiplication on the module $\Pr^k(\E)(U)$ with the property that for all $w$ in $\Pr^k(\E)(U)$
we have 
\[ a(wb)=(aw)b.\]
We write $\Pr^k(\E)^{left}$ (resp. $\Pr^k(\E)^{right}$) to specify we are
considering the left (resp. right) structure.

Note that for $X$ smooth over $S$ and $\E$ locally free
of finite rank, it follows that $\Pr^k_X(\E)$ is locally free of finite rank as left and right $\O_X$-module separately.
%
%
%
%
%
%
\section{Calculation of transition matrices over $\Z$}
In this section we calculate the transition-matrices $L^k_n$ and $R^k_n$ of the principal
parts $\Pr^k(\O(n))$ on the projective line over $\Z$ for all $n\in \Z$ and all $k\geq 1$. 
The structure-matrices defining the principal parts on the projective line are matrices with entries being
Laurent-polynomials with binomial coefficients, hence the principal parts are naturally defined over $\Z$.
The matrices $L^k_n$ and $R^k_n$ are elements of the group $\GL(k+1,\Z[t,t^{-1}])$. The matrix $L^k_n$ define the left $\O$-module 
structure, and the matrix $R^k_n$ define the right $\O$-module structure. Note that in this section all tensor products 
are over $\Z$.

%
%
%
%
%

\textbf{Notation}.
Consider the projective line $\P=\P^1_\Z$ over the integers $\Z$. By definition $\P$ equals $\Proj \Z[x_0,x_1]$. Let 
$U_i=\Spec \Z[x_0,x_1]_{(x_i)}$. Let furthermore
$t=\frac{x_1}{x_0}\otimes 1$ and $u=1\otimes \frac{x_1}{x_0}$ under the isomorphism 
 \[ \Z[\frac{x_1}{x_0}]\otimes_{\Z} \Z[\frac{x_1}{x_0}]\cong \Z[u,t] .\]
It follows from Definition \ref{principalparts}  and Proposition \ref{local} that the sheaf of principal 
parts $\Pr^k(\O(n))|_{U_0}$ as a left $\O_{U_0}$-module equals
\[ \Z[t]\{1\otimes x_0^n, dt\otimes x_0^n, \dots , dt^k\otimes x_0^n\} \]
as a free $\Z[t]$-module where $dt^p=(u-t)^p$. 
Let 
$\t=\frac{x_0}{x_1}\otimes 1$ and $\u=1\otimes \frac{x_0}{x_1}$ under the isomorphism 
 \[ \Z[\frac{x_0}{x_1}]\otimes_{\Z} \Z[\frac{x_0}{x_1}]\cong \Z[\u,\t] .\]
It follows from Definition \ref{principalparts}  and Proposition \ref{local} that  
$\Pr^k(\O(n))|_{U_1}$ as a left $\O_{U_1}$-module equals
\[ \Z[s]\{1\otimes x_1^n,ds\otimes x_1^n, \dots , ds^k\otimes x_1^n\} \]
as a free $\Z[s]$-module, where $s=\frac{1}{t}$ and $ds^p=(\u-\t)^p$. 
Let 
\[ C=\{1\otimes x_0^n,dt\otimes x_0^n, \dots , dt^k\otimes x_0^n\} \]
and 
\[ C'=\{1\otimes x_1^n,ds\otimes x_1^n, \dots , ds^k\otimes x_1^n\} .\]
We aim to describe the principal parts by calculating the transition-matrix 
\[ L^k_n=[I]^{C'}_C\]
between
the bases $C'$ and $C$ as an element of the group $\GL(k+1,\Z[t,t^{-1}])$. 
The following theorem gives the structure-matrix $L^k_n$ describing the left $\O$-module structure
of the principal parts on the projective line over $\Z$. 
\begin{theorem} \label{leftstructure} Consider $\Pr^k(\O(n))$ as left $\O$-module on $\P$. 
In the case $n<0$ the structure-matrix $L^k_n$ is given by the following formula:
\begin{equation}\label{left1}
ds^p\otimes x_1^n=\sum_{j=p}^{k}(-1)^{j}\binom{j-n-1}{p-n-1}t^{n-p-j} dt^{j}\otimes x_0^n ,
\end{equation}
where $0\leq p \leq k$. If $n\geq 0$ the matrix $L^k_n$ is given as follows:

Case $1\leq k \leq n$: 
\begin{equation} \label{left2}
ds^p\otimes x_1^n= \sum_{j=p}^{k}(-1)^p \binom{n-p}{j-p}t^{n-p-j}dt^{j}\otimes x_0^n, 
\end{equation}
where $0\leq p \leq k$. 

Case $n<k$: 
\begin{equation}\label{left3}
ds^p\otimes x_1^n=\sum_{j=p}^{n}(-1)^p\binom{n-p}{j-p}t^{n-p-j}dt^{j}\otimes x_0^n, 
\end{equation}
if $1\leq p \leq n$, and 
\begin{equation}\label{left4}
ds^p\otimes x_1^n=\sum_{j=p}^{k}(-1)^{j}\binom{j-n-1}{p-n-1}t^{n-p-j}dt^{j}\otimes x_0^n, 
\end{equation}
if $n<p\leq k$. 
\end{theorem} 
\begin{proof} Case \ref{left1}: let $0\leq p \leq k$.
\[ ds^p\otimes x_1^n=(\u-\t)^pu^n\otimes x_0^n=(-1)^p\frac{1}{t^p}dt^p(t+dt)^{n-p}\otimes x_0^n=\]
\[ (-1)^p\frac{1}{t^p}dt^p\frac{1}{(t+dt)^{p-n}}\otimes x_0^n =
(-1)^p\frac{1}{t^p}dt^p\frac{1}{t^{p-n}(1+dt/t)^{p-n}}\otimes x_0^n =\]

We have the formula
\[ \frac{1}{(1+x)^l}=\sum_{i\geq 0}(-1)^i\binom{i+l-1}{l-1}x^i \]
and if we let $x=dt/t$ and $l=p-n$ we get
\[
(-1)^p\frac{1}{t^{2p-n}}dt^p\sum_{i\geq0}(-1)^i\binom{i+p-n-1}{p-n-1}\frac{1}{t^i}dt^i\otimes
x_0^n =\]
\[ \sum_{i=0}^{k-p}(-1)^{i+p}\binom{i+p-n-1}{p-n-1}t^{n-2p-i}dt^{i+p}\otimes x_0^n.\]
Changing index with $j=i+p$ we get
\[ \sum_{j=p}^k(-1)^j\binom{j-n-1}{p-n-1}t^{n-p-j}dt^j\otimes x_0^n, \]
and equation \ref{left1} is proved. We leave the proofs of equations
\ref{left2}-\ref{left4} as an excercise.

\end{proof} 
%
%
%
We use similar calculations as in the Theorem \ref{leftstructure} to study the structure-matrix $R^k_n$ defining the principal
parts as right $\O$-module on the projective line over $\Z$.

\textbf{Notation}. It follows from Definition \ref{principalparts} and Proposition \ref{local} that the 
principal parts $\Pr^k(\O(n))|_{U_0}$ as a right $\O_{U_0}$-module equals
\[ \Z[u]\{1\otimes x_0^n,du\otimes x_0^n, \dots , du^k\otimes x_0^n\} \]
as a free $\Z[u]$-module, where $du^p=(t-u)^p$. 
It follows that  $\Pr^k(\O(n))|_{U_1}$ as a right $\O_{U_1}$-module equals
\[ \Z[v]\{1\otimes x_1^n,dv\otimes x_1^n, \dots , dv^k\otimes x_1^n\} \]
as a free $\Z[v]$-module, where $v=\frac{1}{u}$ and $dv^p=(\t-\u)^p$. 
Let 
\[ D=\{1\otimes x_0^n,du\otimes x_0^n, \dots , du^k\otimes x_0^n\} \]
and 
\[ D'=\{1\otimes x_1^n,dv\otimes x_1^n, \dots , dv^k\otimes x_1^n\} .\]
We wish to calculate the transition-matrix 
\[ R^k_n=[I]^{D'}_D\]
describing the right module structure on the principal parts. 
The following theorem gives the structure-matrix $R^k_n$:
\begin{theorem} \label{rightstructure} Consider $\Pr^k(\O(n))$ as right $\O$-module on $\P$. If $n\in \Z$ and $k\geq 1$ 
the following formulas hold: 
\begin{equation} \label{right1}
 1\otimes x_1^n=u^n\otimes x_0^n 
\end{equation} 
and
\begin{equation}\label{right2}
 dv^p\otimes x_1^n=\sum_{j=p}^{k}(-1)^{j}\binom{j-1}{p-1}u^{n-p-j}du^{j}\otimes x_0^n, 
\end{equation}
where $1\leq p \leq k$. 
\end{theorem} 
\begin{proof} Equation \ref{right1} is trivial hence we prove equation \ref{right2}: Let $1\leq p \leq k$. 
We get
\[ dv^p\otimes x_1^n= (\t-\u)^pu^n\otimes x_0^n=(-1)^pu^{n-p}du^p\frac{1}{(u+du)^p}\otimes x_0^n =\]
\[ (-1)^pu^{n-2p}du^p\frac{1}{(1+\u du)^p}\otimes x_0^n .\]
Since we have the formula 
\[ \frac{1}{(1+x)^l}=\sum_{i\geq 0}(-1)^i\binom{i+l-1}{l-1}x^i \]
We get
\[ (-1)^pu^{n-2p}du^p\sum_{i=0}^k
(-1)^i\binom{i+p-1}{p-1}\frac{1}{u^i}du^i\otimes x_0^n .\]
Multiplying we get
\[ \sum_{i=0}^{k-p}(-1)^{i+p}\binom{i+p-1}{p-1}u^{n-2p-i}du^{i+p}\otimes x_0^n. \]
Letting $j=i+p$ we get
\[  \sum_{j=p}^{k}(-1)^{j}\binom{j-1}{p-1}u^{n-p-j}du^{j}\otimes x_0^n, \]
and equation \ref{right2} is proved, and the theorem follows.
\end{proof} 
%
%
%
%
%
The Theorems \ref{leftstructure} and \ref{rightstructure} describe completely the matrices $L^k_n$ and $R^k_n$ 
for all $n\in \Z$ and $k\geq 1$. By \cite{GRO} Theorem 2.1 and \cite{harder} Theorem 3.1 we know that any finite rank locally free
sheaf on $\P^1$ over any field splits into a direct sum of invertible sheaves. The formation of principal parts
commutes with direct sums, hence it follows from Theorems \ref{leftstructure} and \ref{rightstructure} that we have
calculated the transition matrix of $\Pr^k(\E)$ for any locally free finite rank sheaf $\E$ as left and right
$\O$-module on the projective line over any field.

\section{Maps of sheaves and linear equations} 
In this section we develop criteria to split the principal parts on
the projective line over arbitrary rings. We prove existence of
systems of linear equations with the property that solutions to these
systems with coefficients in a ring $A$ gives rise to a splitting of
the principal parts on $\p^1_A$ - the projective line over $A$. Let
$\p^1_A=\Proj(A[x_0,x_1])$
and let $t=x_1/x_0$ and $s=1/t$. 
It follows that
\[ \Pr^k(\O(n))|_{U_0}=A[t]\{1\otimes x_0^n,\cdots, dt^k\otimes x_0^n
\} \]
and
\[ \Pr^k(\O(n))|_{U_1}=A[s]\{1\otimes x_1^n,\cdots ,ds^k\otimes x_1^n
\} ,\] 
 where $U_i\subseteq \p^1_A$ with $i=0,1$ are the basic open subsete of the
 projective line over $A$.
Let in the following $n<0$, $k\geq 1$ and recall the formula from
 Theorem \ref{leftstructure}(with $0\leq p \leq k$):
\begin{align} \label{formula} ds^p\otimes x_1^n=\sum_{j=p}^k
(-1)^j\binom{j-n-1}{p-n-1}t^{n-p-j}dt^j\otimes x_0^n .
\end{align} 
Let furthermore for $0\leq j \leq i \leq k$ $x^i_j$ be independent
variables over $\Z$ and let $0 \leq l \leq k$. Let also $J=\Pr^k(\O(d))$
Make the following definition:
\[ \phi^1_l(x_1^{n-k})=\sum_{p=0}^l x^l_{l-p}t^{p-l}ds^p\otimes x_1^n
.\]
We seek to define for each $l=0,..,k$ a map of left modules 
\[\phi_l: \O(n-k)\rightarrow J \]
and we define
\[ \phi^0_l(x_0^{n-k})=t^{k-n}x_1^{n-k}= \]
\[ t^{k-n}\sum_{p=0}^l x^l_{l-p}t^{p-l}ds^p\otimes x_1^n =
\sum_{p=0}^l x^l_{l-p}t^{k-n-l+p}ds^p\otimes x_1^n .\]
We seek to give relations between the variables $x^i_j$ in order for
the maps $\phi^0_l$ and $\phi^1_l$ to glue to a well-defined map of sheaves. 
We get using formula \ref{formula} the following: 
\[ \phi^0_k(x_0^{n-k})= \sum_{p=0}^l x^l_{l-p}t^{k-n-l+p}\sum_{j=p}^k
(-1)^j\binom{j-n-1}{p-n-1}t^{n-p-j}dt^j\otimes x_0^n  =  \]
\[ \sum_{p=0}^l\sum_{j=p}^k
(-1)^jx^l_{l-p}\binom{j-n-1}{p-n-1}t^{k-l-j}dt^j\otimes x_0^n. \]
By convention
\[ \binom{j-n-1}{p-n-1}=0 \]
when $j<p$ hence we get the expression
\[ \sum_{p=0}^l\sum_{j=0}^k
(-1)^jx^l_{l-p}\binom{j-n-1}{p-n-1}t^{k-l-j}dt^j\otimes x_0^n = \]
\[ \sum_{j=0}^k\sum_{p=0}^l
(-1)^jx^l_{l-p}\binom{j-n-1}{p-n-1}t^{k-l-j}dt^j\otimes x_0^n. \]
\[ \sum_{j=0}^k(-1)^j(\sum_{p=0}^l
x^l_{l-p}\binom{j-n-1}{p-n-1})t^{k-l-j}dt^j\otimes x_0^n. \]
Define
\[ c^l_j= \sum_{p=0}^l x^l_{l-p}\binom{j-n-1}{p-n-1} .\]

We get the expression
\[ \phi^0_l(x_0^{n-k})=\sum_{j=0}^k(-1)^jc^l_jt^{k-l-j}dt^j\otimes x_0^n = \]
\[\sum_{j=0}^{k-l-1}(-1)^jc^l_jt^{k-l-j}dt^j\otimes
x_0^n+(-1)^{k-l}c^l_{k-l}dt^{k-l}\otimes x_0^n+ \]
\[ (-1)^{k-l+1}c^l_{k-l+1}\frac{1}{t}dt^{k-l+1}\otimes x_0^n +\cdots
+(-1)^kc^l_k\frac{1}{t^k}dt^k\otimes x_0^n. \]
We see that if the system of linear equations 
\begin{align}\label{syst1}
S_{k,n}^l: c^l_{k-l}=1,c^l_{k-l+1}=0,...,c^l_k=0 
\end{align}
holds for all $l=0,..,k$ , we get a well defined map of left modules
\[ \phi_l:\O(n-k)\rightarrow J \] 
defined locally as follows:
\begin{align}\label{map0}
\phi^0_l(x_0^{n-k})=\sum_{j=0}^{k-l-1}(-1)^jc^l_jt^{k-l-j}dt^j\otimes
x_0^n 
\end{align}
and 
\begin{align} \label{map1}
\phi^1_l(x_1^{n-k})=\sum_{p=0}^l x^l_{l-p}t^{p-l}ds^p\otimes x_1^n.
\end{align} 
Note that the system $S_{k,n}^l$ is defined over $\Z$. 
\begin{theorem} \label{splitting_one} Let $n<0$ and $k\geq 1$. If for all 
$l=0,..,k$ the systems $S_{k,n}^l$ from \ref{syst1} has
  solutions $x^i_j$ in a ring $A$ with the property that 
\begin{align} \label{determinant}
 x^0_0x^1_0\cdots x^k_0 \in A^* 
\end{align}
then there exists an isomorphism
\[\phi: \oplus^{k+1}\O(n-k)\rightarrow \Pr^k(\O(n))^{left} \]
as left modules on $\p^1_A$.
\end{theorem}
\begin{proof} Assume the system $S^l_{k,n}$ from \ref{syst1} 
has solutions in $A$ for all $l=0,..,k$. By construction we get a map of left modules
  $\phi=\oplus_{l\geq 0} \phi_l$
\[ \phi:\oplus^{k+1}\O(n-k) \rightarrow \Pr^k(\O(n)) .\]
We want to prove that the map $\phi$ is an isomorphism of left
modules. We first consider the matrix of $\phi^0=\phi|_{U_0}$:
\[\phi^0= \begin{pmatrix} \pm c^0_0t^k       & \pm c^1_kt^{k-1}      &   0   &
  0 & \pm c^k_0 \\
                \pm c^0_1t^{k-1}   & \pm c^1_1t^{k-2}      &   0
                    &         \cdots & 0 \\
                    \vdots         &   \vdots              &   \vdots
                    & \cdots & 0 \\ 
                    \vdots         &   \pm c^1_{k-1}              &   \vdots
                    &  \cdots  & \vdots \\
                \pm  c^0_k         &     0                 &   0
                 & 0    &    0 \end{pmatrix} .\]
We see that by definition $|\phi^0|=1$ which is a unit in $A$, hence the matrix
                defines an isomorphism. 
Consider the matrix $\phi^1=\phi|_{U_1}$:
\[ \phi^1= 
\begin{pmatrix}  x^0_0 & x^1_1\frac{1}{t} & \cdots &
                x^k_k\frac{1}{t^k} \\
                  0 & x^1_0 & \cdots & 0 \\
               \vdots & \vdots & \cdots & \vdots \\
                0 & 0 & \cdots & x^k_0 \\ \end{pmatrix}. 
\]
The determinant $|\phi^1|=x^0_0x^1_0 \cdots x^k_0\in A^*$ hence it is
a unit by definition, hence an isomorphism, and the claim is proved. 
\end{proof} 

We next consider the structure of $J=\Pr^k(\O(n))$ in the case when $0
\leq n <k$. Let $\O$ have local basis $e$ on $U_0$ and $f$ on $U_1$ on
$\P^1_A$ where $A$ is an arbitrary ring. Let furthermore $y^i_j$ be
independent variables over $\Z$ for all $0\leq j\leq i \leq n$. Using
similar calculations as above we may do the following: 
Define for $0\leq l \leq n$
\[ d^l_j=\sum^l_{p=0}(-1)^py^l_{l-p}\binom{n-p}{j-p} .\]
Let furthermore
\[ \phi^0_l(e)=\sum_{j=0}^{n-l}d^l_jt^{n-l-j}dt^j\otimes x_0^n, \]
and
\[ \phi^1_l(f)=\sum_{p=0}^l y^l_{l-p}t^{p-l}ds^p\otimes x_1^n.\]
If the system of equations 
\begin{align}\label{syst2}
S^l_{k,n}: d^l_{n-l}=1,d^l_{n-l+1}=\cdots =d^l_n=0 
\end{align}
is satisfied, we get for each $l=0,..,n$ a map of left modules
\[ \phi_l:\O\rightarrow J, \]
Let furthermore $z^i_j$ be independent variables over $\Z$ for $0 \leq
j \leq i \leq k-n-1$ and define for $0\leq m \leq k-n-1$ 
\[ e^m_j=\sum_{j=0}^m z^m_{m-i}\binom{j-n-1}{i} .\]
We seek to define maps
\[ \psi_m:\O(-k-1)\rightarrow J \]
of left modules. Make the following definitions:
\[ \psi^0_m(x_0^{-k-1})=\sum_{j=n+1}^{k-m}(-1)^je^m_j
t^{k-m-j}dt^j\otimes x_0^n \]
and
\[ \psi^1_m(x_1^{-k-1})=\sum_{i=0}^m z^m_{m-i}t^{i-m}ds^{n+1+i}\otimes
x_1^n .\]
Consider the system of linear equations
\begin{align}\label{syst3}
T^m_{k,n}:c^m_{k-m}=1,c^m_{k-m+1}=\cdots =c^m_k=0 .
\end{align}
\begin{theorem} \label{splitting_two} Let $0\leq n < k$. If the systems $S^l_{k,n}$ and
  $T^m_{k,n}$ from Equations \ref{syst2} and \ref{syst3} have
  solutions in $A$ 
with the property that 
\[ y^0_0y^1_0\cdots y^n_0\in A^* \]
and 
\[ z^0_0z^1_0\cdots z^{k-(n+1)}_0 \in A^* \]
there exists an isomorphism
\[ \psi: \O^{n+1}\oplus \O(-k-1)^{k-n} \cong \Pr^k(\O(n))^{left} \]
as left modules on $\p^1_A$.  
\end{theorem}
\begin{proof} The proof is similar to the proof of Theorem
  \ref{splitting_one} and we leave the details to the reader. 
\end{proof}

We next consider the principal parts as right module. Let $n \in \Z$
and $k\geq 1$. Let $\p^1_A=\Proj(A[x_0,x_1])$ and pick coordinates
$u=\frac{x_1}{x_0}$ on $U_0$ and $v=1/u$ on $U_1$. We have 
isomorphisms
\[ \Pr^k(\O(n))|_{U_0}\cong A[u]\{1\otimes x_0^n,\cdots ,du^k\otimes
x_0^n \} \]
and
\[ \Pr^k(\O(n))|_{U_1}\cong A[v]\{1\otimes x_1^n,\cdots ,dv^k\otimes
x_1^n \} \]
of free modules. 
Recall from Theorem \ref{rightstructure} the formulas
\begin{equation} \label{right1}
 1\otimes x_1^n=u^n\otimes x_0^n 
\end{equation} 
and
\begin{equation}\label{right2}
 dv^p\otimes x_1^n=\sum_{j=p}^{k}(-1)^{j}\binom{j-1}{p-1}u^{n-p-j}du^{j}\otimes x_0^n, 
\end{equation}
where $1\leq p \leq k$. 
Define a map of right modules
\[ \O(n)\rightarrow \Pr^k(\O(n))^{right} \]
as follows:
\[ x_0^n \rightarrow 1\otimes x_0^n \]
and
\[ x^n_1\rightarrow 1\otimes x_1^n .\]
It follows from Equation \ref{right1} that we get a well-defined map
of right modules. We seek to define for $0\leq l \leq k-1$ 
maps 
\[ \psi_l:\O(n-k-1)\rightarrow \Pr^k(\O(n))^{right} \]
of right modules. Let for $0\leq j \leq i \leq k-1$ the variables $w^i_j$ be
independent over $A$. 
Make the folloting definition:
\[ \psi^1_l(x_1^{n-k-1})=\sum_{i=0}^l w^k_{l-i}u^{i-l}dv^{i+1}\otimes
x_1^n . \]
We define
\[ \psi^0_l(x_0^{n-k-1})=x_0^{n-k-1}=u^{k+1-n}x_1^{n-k-1}= \]
\[ u^{k+1-n}\sum_{i=0}^l w^l_{l-i}u^{i-l}dv^{i+1}\otimes x_1^n .\]
Using Equation \ref{right2} we get 
\[
\psi^0_l(x_0^{n-k-1})=u^{k+1-n}\sum_{i=0}^lw^l_{l-i}u^{i-l}\sum_{j=i+1}^k
(-1)^j\binom{j-1}{i} u^{n-i-1-j}du^j\otimes x_0^n =\]
\[ \sum_{i=0}^l\sum_{j=i+1}^k
(-1)^jw^l_{l-i}\binom{j-1}{i}u^{k-l-j}du^j\otimes x_0^n.\]
We let $\binom{a}{b}=0 $ if $a<b$ and get the expression
\[ \sum_{j=1}^k(-1)^j\sum_{i=0}^l
w^l_{l-i}\binom{j-1}{i}u^{k-l-j}du^j\otimes x_0^n. \]
Define
\[ f^l_j=\sum_{i=0}^l w^l_{l-i}\binom{j-1}{i} .\]
We get
\[ \psi^0_l(x_0^{n-k-1})=\sum_{j=1}^k(-1)^jf^l_ju^{k-l-j}du^j\otimes
x_0^n = \]
\[ \sum_{j=1}^{k-l-1}(-1)^jf^l_ju^{k-l-j}du^j\otimes x_0^n
+(-1)^{k-l}f^l_{k-l}du^{k-l}\otimes
x_0^n+\]
\[ (-1)^{k-l+1}f^l_{k-l+1}1/udu^{k-l+1}\otimes x_0^n +\cdots +
(-1)^kf^l_k1/u^kdu^k\otimes x_0^n.\]
Let $0\leq l \leq k-1$ and consider the system of linear equations
\begin{align}\label{syst4}
U^l_{k,n}:f^l_{k-l}=1,f^l_{k-l+1}=\cdots =f^l_k=0 .
\end{align}
If the system $U^l_{k,n}$ from Equation \ref{syst4} 
has a solution in $A$ for all $0\leq l \leq k-1$
we get by definition well-defined maps
\[ \psi_l:\O(n-k-1)\rightarrow \Pr^k(\O(n)) \]
of right modules on $\p^1_A$ hence we get a map
\begin{align} \label{rightmap}
\psi:\O(n)\oplus \O(n-k-1)^k \rightarrow \Pr^k(\O(n)) 
\end{align}
of right modules. 

\begin{theorem} \label{splitting_three} Let $n \in \Z$, $k\geq 1$ and $0 \leq l \leq
  k-1$. If the system $U^l_{k,n}$ from Equation \ref{syst4}
has a solution with coefficients in $A$ satisfying 
\[ w^0_0w^1_0\cdots w^{k-1}_0 \in A^* \]
there is an isomorphism
\[ \psi: \O(n)\oplus \O(n-k-1)^k\rightarrow \Pr^k(\O(n))^{right} \]
of right modules on $\p^1_A$. 
\end{theorem}
\begin{proof} We prove that the map $\psi$ introduced above (see
  \ref{rightmap}) is an
  isomorphism. We first consider the  map $\psi|_{U_0}=\psi^0$:
\[\psi^0= \begin{pmatrix} 1 & 0 & 0 & \cdots & 0 \\
0 & \pm f^0_1u^{k-1} & \pm f^1_1u^{k-2} & \cdots & \pm f^{k-1}_1 \\
0 & \pm f^0_2u^{k-2} & \pm f^1_2u^{k-2} & \cdots & 0 \\ 
\vdots & \vdots & \vdots & \vdots & \vdots \\
0 & \pm f^0_{k-1}u & \pm f^1_{k-1} & \cdots & 0 \\
0 & \pm f^0_k & 0 & \cdots & 0  \end{pmatrix} .\]
It follows that $|\psi^0|=\pm f^0_kf^1_{k-1}\cdots f^{k-1}_1=1$ 
hence $\psi|_{U_0}$ is an isomorphism.
We next consider the map $\psi|_{U_1}=\psi^1$:
\[ \psi^1=
\begin{pmatrix} 1 & 0 & 0 & 0 & \cdots & 0 \\ 
0 & w^0_0 & w^1_11/u & \cdots & \cdots & w^{k-1}_{k-1}\frac{1}{u^{k-1}} \\
0 & 0 & w^1_0 & \cdots & \cdots &  w^{k-1}_{k-2}\frac{1}{u^{k-2}} \\
 \vdots & \vdots & \vdots & \vdots & \vdots & \vdots \\ 
0 & 0 & 0 & \cdots & \cdots & w^{k-1}_0 
\end{pmatrix} .\]
It follows that 
\[ |\psi^1|= w^0_0w^1_0\cdots w^{k-1}_0 \in A^* \]
hence the map is an isomorphism, and the claim follows. 
\end{proof} 

The results of this section generalize the results obtained in
\cite{maa1}, where only the principal parts as left module is
considered.

\section{Application: Splitting over  fields of characteristic zero}
In this section we look at the principal parts $\Pr^k(\O(n))$ as left and right $\O$-module on the projective line
for all $n\in \Z$ and $k\geq 1$ over any field $F$ of characteristic zero. We use the results of section 2 and 
techniques developed in 
section 3 to give the splitting-type of the principal parts as left and right module over any field of 
characteristic zero, settling completely the problem of describing the structure of the principal parts 
on the projective line.
%
%
%
%
We first consider the splitting-type of the principal parts as left module.

First a technical lemma.
Consider the following system of linear equations
$A_l \mathbf{x_l}=b_l$:
\[
\begin{pmatrix} \binom{k+m-l-1}{m-1} & \binom{k+m-l-1}{m} & \binom{k+m-l-1}{m+1} &  \cdots & \binom{k+m-l-1}{m+l-1} \\
\binom{k+m-l}{m-1} & \binom{k+m-l}{m} & \binom{k+m-l}{m+1} &  \cdots & \binom{k+m-l}{m+l-1} \\
 \vdots & \vdots & \vdots & \vdots \\
\binom{k+m-2}{m-1} & \binom{k+m-2}{m} & \binom{k+m-2}{m+1} &  \cdots 
& \binom{k+m-2}{m+l-1} \\
\binom{k+m-1}{m-1} & \binom{k+m-1}{m} & \binom{k+m-1}{m+1} &  \cdots & \binom{k+m-1}{m+l-1} 
\end{pmatrix}
\begin{pmatrix} x^l_l \\ x^l_{l-1}  \\ \vdots \\ x^l_{1}  \\ x^l_0  \end{pmatrix}
= \begin{pmatrix} 1 \\ 0 \\ \vdots \\ 0 \\ 0 \end{pmatrix}   ,\]
with $l=0,\dots k$. 
\begin{lemma} \label{lemma} If the vector $\mathbf{x_l}$ is a solution to the
  system defined above, it follows that $x^l_0\neq 0$. 
\end{lemma} 
\begin{proof} Assume on the contrary that $x^l_0=0$. We get the
  following system of linear equations:
\[
\begin{pmatrix} 
\binom{k+n-l}{n-1} & \binom{k+n-l}{n} & \binom{k+n-l}{n+1} &  \cdots & \binom{k+n-l}{n-1+l-1} \\
 \vdots & \vdots & \vdots & \vdots \\
\binom{k+n-2}{n-1} & \binom{k+n-2}{n} & \binom{k+n-2}{n+1} &  \cdots 
& \binom{k+n-2}{n-1+l-1} \\
\binom{k+n-1}{n-1} & \binom{k+n-1}{n} & \binom{k+n-1}{n+1} &  \cdots & \binom{k+n-1}{n-1+l-1} 
\end{pmatrix}
\begin{pmatrix} x^l_l \\ x^l_{l-1}  \\ \vdots \\ x^l_{1}   \end{pmatrix}
= \begin{pmatrix} 0 \\ 0 \\ \vdots \\ 0  \end{pmatrix}   ,\]

\[  \mathbf{Ax}=\mathbf{b} .\]

Using Lemma \ref{det1} with $a=k+n-l, b=n-1$ and $l-1$ we get 
\[ |A|=\frac{ \prod_{i=0}^{l-1} \binom{k+n-l+i}{n-1+i} }{
  \prod_{i=0}^l \binom{k-l+1+i}{i} }\neq 0 \]
and the lemma is proved. 
\end{proof}

\begin{theorem} \label{leftspl1} Consider the principal parts $\Pr^k(\O(n))$ as left module on $\P^1_F$
with $n<0$ and $k\geq 1$ or $1 \leq  k \leq n$. 
There exists an isomorphism 
\[ \phi:\O(n-k)^{k+1}\rightarrow \Pr^k(\O(n))^{left} \]
of left $\O_{\P^1_F}$-modules.
\end{theorem} 
\begin{proof} The case where $1\leq k \leq  n$ is proved in \cite{maa1} Proposition 6.3. 
Consider the case where $n<0$ and $k\geq 1$. Let $0\leq l \leq
k$. According to Theorem \ref{splitting_one}, there exists for all $0\leq l
\leq k$ systems $S^l_{k,n}$ (see \ref{syst1}) of linear equations with the property
that solutions in $F$ for all $0\leq l \leq k$ gives rise to a map
\[\phi :\O(n-k)^{k+1}\rightarrow \Pr^k(\O(n)) \]
of left $\O_{\P^1_F}$-modules. 
We seek solutions in $F$ to the following system of equations
$S^l_{k,n}: A_l\mathbf{x_l}=b_l$:
\[
\begin{pmatrix} \binom{k+m-l-1}{m-1} & \binom{k+m-l-1}{m} & \binom{k+m-l-1}{m+1} &  \cdots & \binom{k+m-l-1}{m+l-1} \\
\binom{k+m-l}{m-1} & \binom{k+m-l}{m} & \binom{k+m-l}{m+1} &  \cdots & \binom{k+m-l}{m+l-1} \\
 \vdots & \vdots & \vdots & \vdots \\
\binom{k+m-2}{m-1} & \binom{k+m-2}{m} & \binom{k+m-2}{m+1} &  \cdots 
& \binom{k+m-2}{m+l-1} \\
\binom{k+m-1}{m-1} & \binom{k+m-1}{m} & \binom{k+m-1}{m+1} &  \cdots & \binom{k+m-1}{m+l-1} 
\end{pmatrix}
\begin{pmatrix} x^l_l \\ x^l_{l-1}  \\ \vdots \\ x^l_{1}  \\ x^l_0  \end{pmatrix}
= \begin{pmatrix} 1 \\ 0 \\ \vdots \\ 0 \\ 0 \end{pmatrix}   ,\]
for all $l=0,\dots k$. 
By Lemma \ref{det1} the determinant of $A_l$ is 
\[ |A_l^1|=\frac{\binom{k-l+m-1}{m-1}\binom{k-l+m}{m}\cdots \binom{k+m-1}{m-1+l}}{\binom{k-l}{0}\binom{k-l+1}{1}\cdots 
\binom{k}{l} } \]
which is non-zero for all $l$ since the characteristic of $F$ is zero, hence there exists  a map 
\[ \phi:\O(n-k)^{k+1}\rightarrow \Pr^k(\O(n)) \]
of left $\O$-modules. 
Consider the maps $\phi^i=\phi|_{U_i}$ for
$i=0,1$. By the proof of Theorem \ref{splitting_one} the determinant of $\phi^0$ equals 
\[ |\phi^0|=\pm c^0_kc^1_{k-1}c^2_{k-2}\cdots c^k_0  \]
which is nonzero by construction, hence $\phi^0$ is an isomorphism.
By the proof of Theorem \ref{splitting_one} the determinant $|\phi^1|$ equals
\[ |\phi^1|= x^0_0x^1_0x^2_0\cdots x^k_0  \]
which is nonzero by Lemma \ref{lemma}. It follows that $\phi^1$ is an isomorphism, and the theorem follows.
\end{proof}

\begin{theorem} \label{leftspl2} Consider the principal parts $\Pr^k(\O(n))$ as left
$\O_{\P^1_F}$-module on $\P^1_F$ with $0 \leq n < k$. There exists an isomorphism 
\begin{equation}\label{iso2}
\phi:\O^{n+1} \oplus \O(-k-1)^{k-n} \rightarrow  \Pr^k(\O(n))^{left} 
\end{equation}
of left $\O_{\P^1_F}$-modules.
\end{theorem}
\begin{proof} By Theorem \ref{splitting_two} there exists systems
  $S^l_{k,n}$ and $T^m_{k,n}$ (see \ref{syst2} and \ref{syst3})  with
  the property that solutions in $F$ gives rise to a map 
\[ \phi: \O^{n+1}\oplus \O(-k-1) \rightarrow \Pr^k(\O(n)) \]
of left modules. 
We seek solutions in $F$ to the following systems of linear equations:
$S^l_{k,n}: A_l\mb{x_l}=b_l$
\[
\begin{pmatrix} \binom{n}{n} & -\binom{n-1}{n-1} & \binom{n-2}{n-2} &  \cdots & (-1)^l\binom{n-l}{n-l} \\
\binom{n}{n-1} & -\binom{n-1}{n-2} & \binom{n-2}{n-3} &  \cdots & (-1)^l\binom{n-l}{n-l-1} \\
 \vdots & \vdots & \vdots & \vdots \\
\binom{n}{n-l+1} & -\binom{n-1}{n-l} & \binom{n-2}{n-l-1} &  \cdots  & (-1)^l\binom{n-l}{n-2l+1} \\
\binom{n}{n-l} & -\binom{n-1}{n-l-1} & \binom{n-2}{n-l-2} &  \cdots & (-1)^l\binom{n-l}{n-2l} 
\end{pmatrix}\begin{pmatrix} y^{l}_{l} \\  y^{l}_{l-1} \\ \vdots \\  y^{l}_1 \\  y^{l}_0 \end{pmatrix}=
\begin{pmatrix} 1 \\ 0 \\ \vdots \\ 0 \\ 0 \end{pmatrix} 
\]
and the system $T^m_{k,n}: B_m\mb{y_m}=c_m$ 
\[
\begin{pmatrix} \binom{k-n-1}{0} & \binom{k+n-1}{1} & \binom{k-n-1}{2} &  \cdots & \binom{k-n-1}{m} \\
\binom{k-n-2}{0} & \binom{k-n-2}{1} & \binom{k-n-2}{2} &  \cdots & \binom{k-n-2}{m} \\
 \vdots & \vdots & \vdots & \vdots & \vdots \\
\binom{k-n-m}{0} & \binom{k-n-m}{1} & \binom{k-n-m}{2} &  \cdots  & \binom{k-n-m}{m} \\
\binom{k-n-m-1}{0} & \binom{k-n-m-1}{1} & \binom{k-n-m-1}{2} &  \cdots & \binom{k-n-m-1}{m} 
\end{pmatrix}
\begin{pmatrix} z^{m}_m \\ z^{m}_{m-1} \\ \vdots \\ z^{m}_1 \\ z^{m}_0 \end{pmatrix}
=\begin{pmatrix} 1 \\ 0 \\ \vdots \\ 0 \\ 0 \end{pmatrix}
,\]
where $l=0,\dots ,n$ and $m=0,\dots ,k-n-1$. The determinant of $A_l$ is by Lemma \ref{det2} equal to $(-1)^{l+1}$,
hence the system $A_l\mb{x_l}=b_l$ has solutions in any field $F$ of characteristic zero.
The determinant of the matrix $B_m$ equals
\[ |B_m|= \frac{ \binom{k-n-1}{0}\binom{k-n}{1}\cdots \binom{k-n+m-1}{m} }{\binom{k-n-1}{0}\binom{k-n}{1}\cdots 
\binom{k-n+m-1}{m} }=1 \]
which is different from zero in $F$. It follows that the system $B_m\mb{y_m}=c_m$ always has solutions in $F$, 
hence we get a map 
\[ \phi: \O^{n+1}\oplus \O(-k-1)^{k-n}\rightarrow \Pr^k(\O(n)) \]
of left $\O_{\P^1_F}$-modules defined as follows: Following the
discussion preceeding Theroem  \ref{splitting_two} and
we define maps $\phi^i_l$  for $l=0,\dots , n$ on $U_0$ as follows:
\begin{equation} \label{map120}
\phi^{0}_l(e)=\sum_{j=0}^{n-l} d^l_jt^{n-l-j}dt^j\otimes x_0^n 
\end{equation}
On $U_1$ we define:
\begin{equation}\label{map121}
\phi^{1}_l(f)=\sum_{p=0}^l y^l_{l-p} t^{p-l}ds^p\otimes x_1^n. 
\end{equation}
We define for $0\leq m \leq k-n-1$ maps $\psi^i_l$ as follows:
\begin{equation} \label{map130}
\psi^{0}_m(x_0^{-k-1})=\sum_{j=k+1}^{k-m}(-1)^j e^m_j t^{k-m-j}dt^j\otimes x_0^n
\end{equation}
On $U_1$:
\begin{equation}\label{map131}
\psi^{1}_l(x_1^{-k-1})=\sum_{i=0}^m z^m_{m-i}  t^{i-m}ds^{n+1+i}\otimes x_1^n.
\end{equation}

Let $\phi=\oplus_{i=0}^{n} \phi_l\oplus_{i=0}^{k-n-1}\psi_i$. We claim that $\phi$ is an isomorphism. 
We consider the map $\phi|_{U_0}$: Its determinant equals
\[ |\phi|_{U_0}|=\pm d^{0}_n d^{1}_{n-1}\cdots d^{n}_{0} e^{0}_{k}e^{1}_{k-1}\cdots e^{k-n-1}_{n+1} \]
which is different from zero by construction, hence $\phi|_{U_0}$ is an isomorphism. 
We next consider the map $\phi_{U_1}$. Its determinant equals
\[ |\phi|_{U_1}|= y^{0}_0 y^{1}_0\cdots  y^{n}_0 z^{0}_0 z^{1}_0 z^{k-n-1}_0 \]
which is nonzero by an argument similar to the one in Lemma
\ref{lemma}, 
hence $\phi|_{U_1}$ is an isomorphism, and the theorem is proved. 
\end{proof}

We get a similar theorem on the splitting-type of the principal parts as right module:

\begin{theorem} \label{rightspl1} Let $F$ be a field of characteristic zero. Consider $\Pr^k(\O(n))$ as right $\O$-module
on $\P^1_F$ for all $n\in \Z$ and $k\geq 1$. There exists an isomorphism
\[ \psi: \O(n)\oplus \O(n-k-1)^k \rightarrow \Pr^k(\O(n))^{right} \]
of right $\O$-modules.
\end{theorem} 
\begin{proof} Define a map
\[ \psi_{-1}:\O(n)\rightarrow \Pr^k(\O(n))^{right} \]
as follows:
\[ \psi_{-1}(x_0^n)=1\otimes x_0^n \]
and
\[ \psi_{-1}(x_1^n)=1\otimes x_1^n \]
According to Theorem \ref{splitting_three} we need to
  prove existence of solutions of the systems $U^l_{k,n}$ (see
  \ref{syst4}) with coefficients in $F$ for all $0 \leq l \leq k-1$.

Existence of the map $f$ is by definition equivalent to the existence of solutions in 
$F$ of the following systems of linear equations $A_l\mb{x_l}=b_l:$ 
\[
\begin{pmatrix} \binom{k-1}{0} & \binom{k-1}{1} & \binom{k-1}{2} &  \cdots & \binom{k-1}{l} \\
\binom{k-2}{0} & \binom{k-2}{1} & \binom{k-2}{2} &  \cdots & \binom{k-2}{l} \\
 \vdots & \vdots & \vdots & \vdots & \vdots \\
\binom{k-l}{0} & \binom{k-l}{1} & \binom{k-l}{2} &  \cdots  & \binom{k-l}{l} \\
\binom{k-l-1}{0} & \binom{k-l-1}{1} & \binom{k-l-1}{2} &  \cdots & \binom{k-l-1}{l} 
\end{pmatrix} \begin{pmatrix} w^{l}_l \\  w^{l}_{l-1}  \\ \vdots \\  w^{l}_1 \\  w^{l}_0 \end{pmatrix}
\begin{pmatrix} 1 \\ 0 \\ \vdots \\ 0 \\ 0 \end{pmatrix} 
,\]
where $l=0,\dots ,k-1$. 
The determinant of the matrix $A_l$ equals
\[ |A_l|= \frac{ \binom{k-1}{0}\binom{k}{1} \cdots \binom{k+l-1}{l} }{\binom{k-1}{0}\binom{k}{1} \cdots \binom{k+l-1}{l} } =1\]
which is a unit, hence the map $\phi$ exists over the field $F$. 

According to the discussion preceeding Theorem \ref{splitting_three} we may
define for $0 \leq l \leq k-1$ maps $\psi_l$ as follows: On $U_0$ we define
\begin{equation}\label{map140}
\psi^0_l(x_0^{n-k-1})=\sum_{j=0}^{k-l}(-1)^j f^l_j u^{k-l-j}du^j\otimes x_0^n.
\end{equation}
Define a map on $U_1$: 
\begin{equation}\label{map141}
\psi^1_l(x_1^{n-k-1})=\sum_{i=0}^l w^l_{l-i} u^{i-l}dv^{i+1}\otimes x_1^n,
\end{equation}
Define $\psi=\oplus_{i=-1}^{k-1}\psi_i$. The determinant of $\phi|_{U_0}$ equals 
\[ |\psi|_{U_0}|=\pm  f^0_kf^1_{k-1}\cdots f^k_0  \]
which is non-zero by construction, hence $\psi|_{U_0}$ is an isomorphism.
The determinant of $\psi|_{U_1}$ equals  
\[ |\psi|_{U_1}|= w^0_0w^1_0\cdots w^{k-1}_0 \] 
which is non-zero using an argument similar to the one in Lemma
\ref{lemma}, 
hence $\psi|_{U_1}$ is an isomorphism, and the theorem is proved.
\end{proof} 

Note that the determinant of the matrix $A_l$ from the proof of Theorem \ref{right1} is one, hence 
the map constructed exists over an arbitrary ring. 

The Theorems \ref{leftspl1}, \ref{leftspl2}  and \ref{rightspl1} generalize Proposition 6.3 and Theorem 7.1 in \cite{maa1}, settling
completely the problem of describing the principal parts on the projective line over a field of characteristic
zero. Since formation of principal parts commute with direct sums and because on the projective line over any field 
any finite rank locally free sheaf splits into invertible sheaves (\cite{GRO},\cite{harder}), we have described the splitting-type 
of the principal parts of any locally free finite rank sheaf as left and right module over any field of characteristic 
zero, generalizing results obtained in \cite{maa1},\cite{PIE} and
\cite{dirocco}. The problem of classifying principal parts on the
projective line over an arbitrary field will be adressed in a future paper.

\section{Appendix: calculations of determinants}

In this section we calculate some determinants appearing in section 3 and 4.

\begin{lemma} \label{det1} Assume $a,b$ are integers with $b\leq a$. Then the following is true for all $l\geq 1$:
The determinant of the matrix 
\begin{equation} \label{matrix1}
M=\begin{pmatrix}  \binom{a}{b} &  \binom{a}{b+1} & \binom{a}{b+2}  &  \cdots &   \binom{a}{b+l}    \\
 \binom{a+1}{b}   & \binom{a+1}{b+1}  & \binom{a+1}{b+2}  &  \cdots & \binom{a+1}{b+l}  \\
 \vdots & \vdots & \vdots & \vdots & \vdots \\
 \binom{a+l-1}{b}  & \binom{a+l-1}{b+1}  & \binom{a+l-1}{b+2}  &  \cdots  & \binom{a+l-1}{b+l}  \\
 \binom{a+l}{b}   & \binom{a+l}{b+1}   & \binom{a+l}{b+2}  &  \cdots &  \binom{a+l}{b+l}
\end{pmatrix} 
\end{equation} 
equals 
\begin{equation} \label{formula1} 
|M|=\frac{ \prod_{i=0}^l \binom{a+i}{b+i} }{ \prod_{i=0}^l \binom{a-b+i}{i} }
\end{equation} 
\end{lemma}
\begin{proof} One easily shows the following formula:
\begin{equation} \label{equ}
\binom{a+i}{b+j}-\frac{a+i+1-b-l}{a+i+1}\binom{a+i+1}{b+j}=\frac{l-j}{b+j}\binom{a+i}{b-1+j} .
\end{equation}
Apply equation \ref{equ} to the matrix \ref{matrix1} to get the matrix
\begin{equation} \label{matrix2}
\begin{pmatrix} \frac{l}{b}\binom{a}{b-1} & \frac{l-1}{b+1}\binom{a}{b} & \frac{l-2}{b+2}\binom{a}{b+1}  &  
\cdots &   \frac{1}{b+l-1}\binom{a}{b+l-2} & 0    \\
 \frac{l}{b}\binom{a+1}{b-1}   &\frac{l-1}{b+1} \binom{a+1}{b}  &\frac{l-2}{b+2} \binom{a+1}{b+1}  &  \cdots &
\frac{1}{b+l-1} \binom{a+1}{b+l-2} & 0  \\
 \vdots & \vdots & \vdots & \vdots & \vdots & 0  \\
 \frac{l}{b}\binom{a+l-1}{b-1}  &\frac{l-1}{b+1} \binom{a+l-1}{b}  &\frac{l-2}{b+2} \binom{a+l-1}{b+1}  &  \cdots  &
\frac{1}{b+l-1} \binom{a+l-1}{b+l-2} & 0  \\
 \binom{a+l}{b}   & \binom{a+l}{b+1}   & \binom{a+l}{b+2}  &  \cdots & \binom{a+l-1}{b+l-1} &  \binom{a+l}{b+l} \end{pmatrix}
\end{equation} 
We get the following determinant
\[ \frac{l!}{b(b+1)(b+2)\cdots (b+l-2)}\binom{a+l}{b+l}|N| ,\] 
where $N$ is the matrix
\begin{equation} \label{matrix2}
\begin{pmatrix}  \binom{a}{b-1} &  \binom{a}{b} & \binom{a}{b+1}  &  \cdots &   \binom{a}{b+l-2}    \\
 \binom{a+1}{b-1}   & \binom{a+1}{b}  & \binom{a+1}{b+1}  &  \cdots & \binom{a+1}{b+l-2}  \\
 \vdots & \vdots & \vdots & \vdots & \vdots \\
 \binom{a+l-1}{b-1}  & \binom{a+l-1}{b}  & \binom{a+l-1}{b+1}  &  \cdots  & \binom{a+l-1}{b+l-2}  
\end{pmatrix} 
\end{equation} 
By induction we get the following formula:
\[ \frac{l!}{b(b+1)(b+2)\cdots (b+l-2)}\binom{a+l}{b+l}
\frac{\binom{a}{b-1}\cdots \binom{a+l-1}{b+l-2}}{\binom{a-b+2}{1}\binom{a-b+3}{2}\cdots \binom{a-b+l}{l-1} },\]
and the formula \ref{formula1} follows.
\end{proof} 

\begin{lemma} \label{det2} Let $a$ be an integer with $a\geq 0$. Then the following is true:
The determinant of the matrix 
\begin{equation} \label{matrix3}
P=\begin{pmatrix}  \binom{a}{a} &  \binom{a-1}{a-1} & \binom{a-2}{a-2}  &  \cdots &   \binom{a-l}{a-l}    \\
 \binom{a}{a-1}   & \binom{a-1}{a-2}  & \binom{a-2}{a-3}  &  \cdots & \binom{a-l}{a-l-1}  \\
 \vdots & \vdots & \vdots & \vdots & \vdots \\
 \binom{a}{a-l+1}  & \binom{a-1}{a-l}  & \binom{a-2}{a-l-1}  &  \cdots  & \binom{a-l}{a-2l+1}  \\
 \binom{a}{a-l}   & \binom{a-1}{a-l-1}   & \binom{a-2}{a-l-2}  &  \cdots &  \binom{a-l}{a-2l}
\end{pmatrix} 
\end{equation} 
equals 
\begin{equation} \label{formula2} 
|P|=(-1)^{l+1}.
\end{equation} 
\end{lemma}
\begin{proof} The matrix $P$ equals the matrix 
\begin{equation} \label{matrix4}
P=\begin{pmatrix}  1  &   1     &  1       &  \cdots &    1       \\
  a   &   a-1   &  a-2   &  \cdots &  a-l  \\
  \frac{a(a-1)}{2!}  & \frac{(a-1)(a-2)}{2!}  & \frac{(a-2)(a-3)}{2!}  & \cdots & \frac{(a-l)(a-l-1)}{2!} \\
 \vdots & \vdots & \vdots & \vdots & \vdots \\
 
 \binom{a}{a-l}   & \binom{a-1}{a-l-1}   & \binom{a-2}{a-l-2}  &  \cdots &  \binom{a-l}{a-2l}
\end{pmatrix} .
\end{equation} 
We get the Vandermonde-matrix
\begin{equation} \label{matrix5}
\frac{1}{1!2!\cdots l!}\begin{pmatrix}  1  &   1     &  1       &  \cdots &    1       \\
  a   &   a-1   &  a-2   &  \cdots &  a-l  \\
  a^2  & (a-1)^2  & (a-2)^2  & \cdots & (a-l)^2 \\
 \vdots & \vdots & \vdots & \vdots & \vdots \\
 
 a^l & (a-1)^l & (a-2)^l & \cdots & (a-l)^l 
\end{pmatrix},
\end{equation} 
which is easily seen to have determinant $(-1)^{l+1}$, and the lemma follows.
\end{proof}

\textbf{Acknowledgments}. This paper was written in the period autumn 2003-
spring 2004 and I want to thank Dan Laksov for comments on 
the problems discussed in this paper. Thanks also to Andrei Zelevinsky.
It is also a great pleasure to thank colleagues and staff at the Royal Institute of Technology for providing a 
friendly and inspiring atmosphere during its preparation.

\end{document}